\documentclass[]{iopart}
\usepackage{iopams}
\usepackage{setstack}
\usepackage[]{graphicx}
\usepackage{color}
\usepackage[latin1]{inputenc}

\newcommand{\colr}{}

\newcommand{\Prob}{{\rm Prob}}
\newcommand{\const}{{\rm const.}}

\newcommand{\au}{{\underline{a}}}
\newcommand{\ao}{{\overline{a}}}

\begin{document}

\title[]{Interval estimations in metrology}
\author{G Mana$^1$ and C Palmisano$^2$}
\address{$^1$INRIM - Istituto Nazionale di Ricerca Metrologica, Str.\ delle Cacce 91, 10135 Torino, Italy}
\address{$^2$UNITO - Università di Torino, Dipartimento di Fisica, v. P.\ Giuria, 1 10125 Torino, Italy}
\ead{g.mana@inrim.it}

\begin{abstract}
This paper investigates interval estimation for a measurand that is known to be positive. Both the Neyman and Bayesian procedures are considered and the difference between the two, not always perceived, is discussed in detail. A solution is proposed to a paradox originated by the frequentist assessment of the long-run success rate of Bayesian intervals.
\end{abstract}
\submitto{Metrologia}
\pacs{02.50.Cw, 02.50.Tt, 06.20.Dk, 07.05.Kf}
% 02.50.Cw Probability theory
% 02.50.Tt Inference methods
% 06.20.Dk Measurement and error theory
% 07.05.Kf Data analysis: algorithms and implementation; data management

%\baselineskip 8mm

\section{Introduction}
{\colr The Neyman and Bayes viewpoints about how to carry out interval estimation \cite{Neyman:1935,Neyman:1937,Jaynes:1976,Sivia} lead to different uncertainty statements. Since they are calculating different intervals, there is a debate over the meaning of confidence level and coverage probability for an uncertainty statement.} A non-exhaustive list of papers investigating this issue is \cite{Stein:1959,Feldman:1998,CERN:2000,Bukin:2003,Lira:2006,Hall:2008,Lira:2008,Wang-Iyer:2009,Possolo:2009,Willink:2010,Willink:2010c,Attivissimo:2012,Mana:2013}. A simple problem that makes the viewpoints' differences evident is when there is a measurement of a real quantity that is small with respect to the measurement uncertainty and that it is known to be have a well defined sign \cite{Calonico:2008,Calonico:2009}. In \cite{Willink:2010a,Willink:2010b} Willink showed that, in a Monte Carlo simulation of repeated Gaussian measurements of the same measurand, the long-run success rate of Bayesian intervals to encompass the measurand disagrees with the expected value. This paper presents the results of an investigation designed to understand the basic concepts of interval estimation and to explain this paradoxical result.

When reporting the uncertainty of measurements, the awareness of the differences between the Neyman and Bayesian approaches is essential. Interval estimation is a procedure to find a pair of values that succeeds in including the measurand with a stipulated probability. In the Neyman approach, the focus is on different interval estimations, given the same measurand value. In the Bayesian approach, the focus is on different measurand values, given the same interval. After reviewing the Neyman and Bayesian solutions, the paper illustrates interval estimation given a Gaussian sample of a positive quantity. {\colr Eventually, it shows that both approaches achieve the stipulated success rate.}

\section{Interval estimation}
Before the measurement is carried out, the measured value, $x$, can be viewed as an unknown member of a population described by a probability distribution, $P_x(\xi|a)$, parameterised by the measurand value $a$, which -- though unknown -- has a fixed value. In $P_x(\xi|a)$, the letter $\xi$ is a dummy variable which labels the space of the possible $x$ values; the vertical bar indicates that the probability density in $x=\xi$ is conditioned on a measurand value equal to $a$.

The probability distribution of the measurand values enters our considerations because we face a range of possible values, but we are not able to figure out what it is. Probability assignments to the $a$ values and probability calculus make our knowledge quantitative and allow us to come to sensible decisions. Prior to the measurement, the measurand value can be viewed as an unknown member of a population described by the probability distribution $\pi(\phi)$ -- where $\phi$ labels the possible $a$ values. When the measured value $x_0$ is on hand, $\pi(\phi)$ must be updated to $P_a(\phi|x_0)$, where $x_0$ is a known parameter. These distributions are linked by the Bayes theorem
\begin{equation}\label{bayes}
 P_a(\phi|x_0) = \frac{L(\phi;x_0)\pi(\phi)}{Z(x_0)} ,
\end{equation}
where $L(\phi;x_0)=P_x(x_0|\phi)$ is the likelihood function and $Z(x_0)$ is a normalising factor.

Given the measurement result $x_0$, interval estimation is the problem of finding an interval -- $[a_1,a_2]$, which is called a credible (or coverage) interval -- such that the measurand value in it with a predetermined probability -- $\Prob(a\in[a_1,a_2]|x_0)$, which is called coverage probability and is represented by $\alpha$. Therefore, the interval end-points are the solutions of
\begin{equation}\label{Pax0}
 \Prob(a\in[a_1,a_2]|x_0) = \alpha .
\end{equation}
It must be noted that (\ref{Pax0}) is conditioned on the fixed measured value $x_0$.

\subsection{Neyman: confidence intervals}
According to Neyman, it is meaningless to assign probabilities to the possible measurand value; he discarded (\ref{bayes}) and (\ref{Pax0}) and proposed a statistic -- that is, a function of the measurement result -- having, in a long series of repeated application to different measured values of the same measurand, a success rate of including the measurand equal to the coverage probability \cite{Neyman:1937}. The Neyman interval, $[\au,\ao]$, is called confidence interval and the (predetermined) success rate of the procedure, $\Prob([\au,\ao]\ni a|a)$, is called confidence level. The statement $[\au,\ao]\ni a$ maintains the unknown $[\au,\ao]$ interval includes the known measurand value $a$, whereas $a\in[a_1,a_2]$ maintains the unknown measurand value $a$ is included in the known $[a_1,a_2]$ interval.

In the same way as the measurement result, the confidence interval $[\au,\ao]$ is picked at random from an interval set -- the urn of the frequentist model -- where the fraction of intervals containing the measurand value (the confidence level) is equal to the coverage probability. Therefore, $[\au,\ao]$ is an estimate of the credible interval $[a_1,a_2]$ and the confidence level is a property of the estimator, non of the specific interval sampled. As Neyman repeatedly stated, the confidence level is only the probability that a future interval embeds the measurand.

The interval end-points are the solutions of
\begin{equation}\label{Paa}
 \Prob([\au,\ao]\ni a|a) = \alpha .
\end{equation}
It must be noted that (\ref{Paa}) is conditioned on the fixed measurand value $a$; since it is unknown, (\ref{Paa}) is meaningful only if the statistics used to calculate $\au$ and $\ao$ is such that $\Prob([\au,\ao]\ni a|a)$ is independent of $a$. Afterwards, since (\ref{Paa}) is independent of the actual measurand value, the probability statement about $[\au,\ao]$ can be restated as one about the $a$ value.

The Neyman procedure uses a pair of continuous and monotonic functions of the measurand, $x_1(a)$ and $x_2(a)$, so chosen as $\Prob(x\in [x_1,x_2]|a)=\alpha$. That is,
\begin{equation}\label{CI1}
 F_x(x_2|a) - F_x(x_1|a) = \alpha ,
\end{equation}
where $F_x(\xi|a)$ is the cumulative distribution associated to $P_x(\xi|a)$. Provided the measured value $x_0$ is in the domain of both the inverse functions $x_2^{-1}(\xi)$ and $x_1^{-1}(\xi)$, it follows that $\au=x_2^{-1}(x_0)$ and $\ao=x_1^{-1}(x_0)$ are the sought interval end-points. Equivalently, given the measured value $x_0$, the interval end-points are the solution of
\begin{equation}\label{CI2}
 F_x(x_0|\au) - F_x(x_0|\ao) = \alpha .
\end{equation}

In the same way as a measurand estimate is not the measurand value, in general, the probability of $[\au,\ao]\ni a$ -- where $[\au,\ao]$ is built by solving (\ref{CI2}) and $x_0$ is a given measurement result -- is not equal to the confidence-level \cite{Neyman:1937}; that is, $\Prob([\au,\ao]\ni a|a) = \alpha $ but $\Prob(a\in[\au,\ao]|x_0)\ne\alpha$. Rather, $[\au,\ao] $ is randomly sampled from a set of intervals where the fraction $\alpha$ of them embeds $a$. According to (\ref{Paa}), this sample space is the set of the intervals built by solving (\ref{CI2}), where $x_0$ are the results of repeated measurement of the same measurand.

Some remarks are needed. Firstly, (\ref{CI1}) and (\ref{CI2}) do not identify $[x_1,x_2]$ and $[\au,\ao] $ uniquely; in order to have a single solution, additional constraints are necessary. Secondly, once a constraint has been chosen, there may exist measurement results wherefore (\ref{CI2}) has no solution. Thirdly, when there are nuisance parameters, no general algorithm exists to build a confidence interval for the measurand value only.

\begin{figure}\centering
\includegraphics[width=7cm]{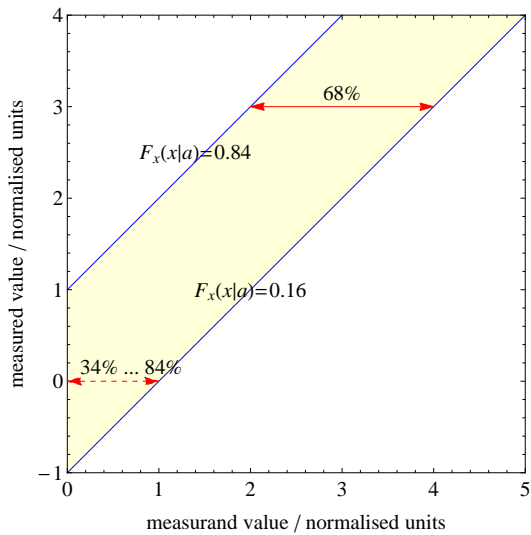}
\caption{16\% and 84\% quantiles of the results of an unbiased Gaussian measurement of a positive quantity. The solid arrow is the solution of (\ref{CI2}) and (\ref{constraint}-$b$), where $x_0/u=3$. The dashed arrow is the solution of (\ref{CI2}) and (\ref{constraint-b}), where the confidence level is any value from 34\% to 84\% and $x_0/u=0$.} \label{upperGauss}
\end{figure}

\subsubsection{Example: Gaussian measurement of a positive quantity.}
The figure \ref{upperGauss} illustrates the Neyman procedure. The measured value $x_0$ of $a>0$ is drawn from the normal distribution $P_x(\xi|a)=N(\xi;a,u^2)$ whose mean and variance are $a$ and $u^2$. For the sake of simplicity, the variance has been set to one, which corresponds to redefine the measurand and measured values as $a/u$ and $x_0/u$. By setting $u^2=1$, the cumulative distribution is
\begin{equation}\label{erfc}
 F_x(\xi|a) = \frac{{\rm erfc}[(a-\xi)/\sqrt{2}]}{2} ,
\end{equation}
where erfc$(x)$ is the complementary error function. In Fig.\ \ref{upperGauss}, the two lines, $x_1=a-u$ and $x_2=a+u$, are so chosen as $\Prob(x\in [x_1,x_2]|a)=0.68$. For the sake of simplicity, the numerical values are rounded to the second digit. Given the measured value $x_0$, the solution of (\ref{CI2}), where
\numparts\begin{eqnarray}\label{constraint}
 F_x(x_0|\au) &= &0.84 \\
 F_x(x_0|\ao) &= &0.16 , \label{constraint-b}
\end{eqnarray}\endnumparts
is $[x_0-u,x_0+u]$; an example is the solid arrow in Fig.\ \ref{upperGauss}.

If $x_0/u < 1$, (\ref{CI2}) has no solution satisfying (\ref{constraint}-$b$). A way to bypass this problem is to allow negative $a$ values and to say that the confidence interval is partly non-physical. However, when $x_0/u < -1$, $[x_0-u,x_0+u]$ is entirely in the negative region and its coverage probability is null; that is, $\Prob(a\in [x_0-u,x_0+u]|x_0)=0$. This is not surprising; besides, a negative interval is not more unusual than a negative datum and $[x_0-u,x_0+u]$ is one of the intervals of the Neyman's sample-space not including $a$. The paradox is solved by observing that it arises only because we know in advance that $a>0$. In addition, it is not correct to identify the procedure confidence-level with the coverage probability of $[x_0-u,x_0+u]$. To say that $[x_0-u,x_0+u]$ is a 68\% confidence interval means that $\Prob([x-u,x+u]\ni a|a)=0.68$, not that $\Prob(a\in [x_0-u,x_0+u]|x_0)=0.68$ \cite{Neyman:1937}. The confidence level is conditioned to the measurand value, not to the measured value. This means that the probability of sampling a future interval such that $[x-u,x+u]\ni a$ is true, where $x$ is unknown and $a$ is known, is 0.68. But, once $x=x_0$ is on hand, the probability -- updated by the information delivered by $x_0$ -- of $a\in [x_0-u,x_0+u]$, where $a$ is unknown and $[x_0-u,x_0+u]$ is known, might be different.

A second solution, proposed in \cite{Willink:2010,Willink:2010a}, is to exclude the $a<0$ values. However, in this case, different confidence levels lead to the same interval. For instance, as shown by dashed arrow in Fig.\ \ref{upperGauss}, given (\ref{constraint-b}), $x_0/u=0$, and any confidence level from 34\% to 84\%, the result is always the $[0,u]$ interval. A further solution is to switch between two-sided intervals and upper limits according to the measured value. For instance, if $x_0/u < 1$, to switch from (\ref{constraint}-$b$) to $\au=0$ and $F_x(x_0|\ao) = 0.32$. However, flip-flopping is inconsistent with a predetermined confidence level. A solution that uses the freedom to choose the $x_1(a)$ and $x_2(a)$ functions is given in \cite{Feldman:1998}. The resulting intervals change continuously from upper limits to two-sided intervals as the measured value becomes more statistically significant.

\subsection{Bayes: credible intervals}\label{bayes-interval}
By definition, the probability of $a\in[a_1,a_2]$ is the integral of $P_a(\phi|x_0)$ between two given limits, $a_1$ and $a_2$. Hence, the end points of credible intervals having a coverage probability equal to $\alpha$ are the solutions of
\begin{equation}\label{CI3}
 F_a(a_2|x_0) - F_a(a_1|x_0) = \alpha ,
\end{equation}
where $F_a(\phi|x_0)$ is the cumulative distribution associated to $P_a(\phi|x_0)$. In the framework of a frequency-of-occurrence model of $\Prob(a\in [a_1,a_2]|x_0)$, the sample space is the set of the $a$ values consistent with the same measurement result and, consequently, with the same credible interval.

It must be noted that, to build credible intervals, the availability of a measurement result is not an essential ingredient. In fact, by resorting to the prior probability distribution $\pi(\phi)$, credible intervals can be built also if no measurement has been carried out. This emphasises again that a probability distribution is not an intrinsic quality of the measurand, but a way to encode our knowledge of its value.

When $\pi(\phi)$ is the uniform distribution and the probability density function of $x$ owns the symmetry $P_x(\xi|\phi)=P_x(\phi|\xi)$, that is, it is invariant with respect to the replacement $\xi \rightleftharpoons \phi$, the Bayes theorem simplifies to $P_a(\phi|x_0) = P_x(\phi|x_0)$ and brings to light that the post-data probability density of the measurand values and the sampling distribution of the measurement results are the same function. If, in addition, $P_x(\xi|\phi)$ is a function of $|\xi-\phi|$ only, it can be proved that the Neyman and Bayesian procedures lead to the same interval. The occurrence of the interval identity -- for instance, when $P_x(\xi|a)$ is the ubiquitous Gaussian distribution -- causes misunderstandings. One may carry out a Neyman interval-estimation and use the result as if $\Prob(a\in[\au,\ao]|x_0)=\alpha$, which, in general, is not correct.

\begin{figure}\centering
\includegraphics[width=7cm]{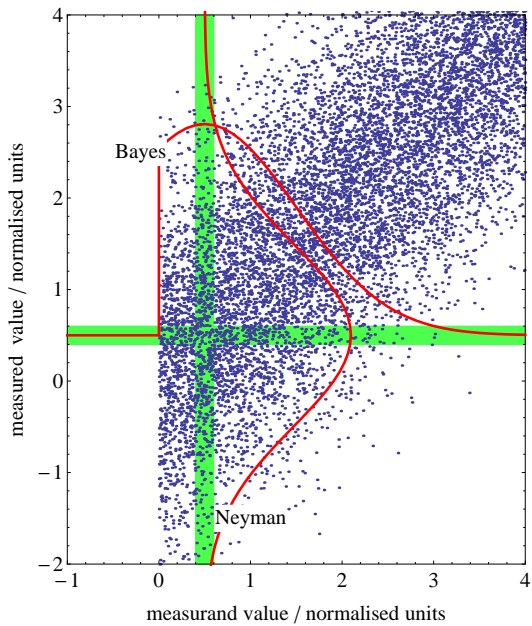}
\caption{Scatter plot of the joint distribution of the $\{a_i,x_i\}$ pairs for a Gaussian measurement of a quantity uniformly distributed in the $[0,\infty]$ interval. The Neyman curve is the sampling distribution $N(\xi;a,u^2)$ of the measurement result for a measurand value $a/u=0.5$. The Bayes curve is the post-data measurand distribution $P_a(\phi|x_0)$, given the measurement result $x_0/u=0.5$. The green strips indicate the sample spaces -- horizontal: $\{a_i, x=\const\}$, fixed measured value; vertical: $\{a=\const,x_i\}$, fixed measurand value -- used to assess the success rates of the Bayes and Neyman solutions of the interval estimation problem.} \label{joint}
\end{figure}

\subsubsection{Example: Gaussian measurement of a positive quantity.}
Let us suppose again that the measured value $x_0$ of $a>0$ is drawn from the normal distribution $P_x(\xi|a)=N(\xi;a,u^2)$. By setting again $u^2=1$, a uniform prior probability distribution of the $a$ values must be updated into
\begin{equation}\label{gauss}
 P_a(\phi|x_0) = \frac{2\exp\big[-(\phi-x_0)^2/2\big]\vartheta(\phi)}{\sqrt{2\pi}\,{\rm erfc}(-x_0/\sqrt{2})} ,
\end{equation}
where erfc$(x)$ is the complementary error function and $\vartheta(\phi)$ is the Heaviside function. The relevant cumulative distribution is
\begin{equation}\label{erf}
 F_a(\phi|x_0) = \frac{{\rm erf}(x_0/\sqrt{2}) + {\rm erf}[(\phi-x_0)/\sqrt{2}]}{{\rm erfc}(-x_0/\sqrt{2})} ,
\end{equation}
where erf$(x)$ is the error function. Given the measured value $x_0$, we will consider the intervals constrained by
\numparts\begin{eqnarray}\label{c-Bayes}
 F_a(a_1|x_0) &= &0.16 \\
 F_a(a_2|x_0) &= &0.84 .
\end{eqnarray}\endnumparts

After a measurement has been completed, the measurement result is a known quantity. Since we own only this unique result, it is not clear what population is to be used to imagine repeated measurements and to build a frequentist model of an unconditioned statement about the probability of the measurand to belong a given interval. In fact, there is a multiplicity of sample spaces to each of which we can regard the unknown $\{a,x\}$ repeated-measurement pairs as belonging, none having an objective reality and all being products of our subjective preference. However, (\ref{Pax0}) and (\ref{Paa}) keep strictly to conditional statements so that the relevant frequentist models can be uniquely defined.

The sample spaces of the frequency-of-occurrence models of $\Prob([\ao,\au]\ni a|a)$ and $\Prob(a\in [a_1,a_2]|x_0)$ are shown in Fig.\ \ref{joint}. The scatter plot shows the joint distribution $P_{x,a}(\xi,\phi)=N(\xi;\phi,u)\pi(\phi)$ of the measurand- and measured-value pairs $\{a_i,x_i\}$ for a Gaussian measurement of a quantity uniformly distributed in the $[0,\infty]$ interval. The pairs $\{a=\const,x_i\}$, having the same measurand value, make up the sample space of $\Prob([\ao,\au]\ni a|a)$ and will be used to assess the success rate of the Neyman intervals: in Fig.\ \ref{joint} they are in the vertical strip. The pairs $\{a_i,x=\const\}$, having the same measured value, make up the sample space of $\Prob(a\in [a_1,a_2]|x_0)$ and will be used to assess the success rate the Bayesian intervals: in Fig.\ \ref{joint} they are in the horizontal strip.

\section{Performance analysis}\label{performance}
This section examines the performances of the Neyman and Bayesian procedures from a frequentist viewpoint. Monte Carlo simulations are used to calculate the success rates of confidence and credible intervals and to compare the results against the expected rates. The case studied is where the measurement of a positive quantity $a$ gives a Gaussian datum having known variance $u^2$.

\subsection{Confidence intervals}
\begin{figure}\centering
\includegraphics[width=7cm]{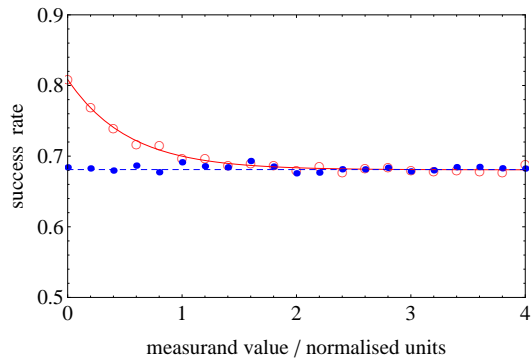}
\caption{Gaussian measurements of a positive quantity: success rate of confidence intervals. Different intervals have been built from the results of repeated measurements of the same measurand. Dots: frequencies observed in Monte Carlo simulations; the horizontal line is the theoretically expected value $\Prob([\au,\ao]\ni a|a)$. The calculations have been carried out by setting (\ref{constraint}-$b$). Empty circles: the intervals entirely in the negative region have been rejected and measurements repeated; the solid line is a smoothed interpolation of the data.} \label{Neyman}
\end{figure}

According to (\ref{Paa}), in a long series of repeated measurements of the same measurand, the fraction $\alpha$ of the different confidence intervals built by solving (\ref{CI2}) contains the single $a$ value. Since it is conditional on a fixed measurand value, to test numerically $\Prob([\au,\ao]\ni a|a)=\alpha$, the measurand value (say, $a_0>0$) must be fixed. Next, a number of measurement results are repeatedly sampled according to $P_x(\xi|a_0)$ and the relevant confidence intervals are built. Each trial involves determining if the interval contains the fixed measurand value. As shown by the horizontal line in Fig.\ \ref{Neyman}, a Monte Carlo simulation, carried out by fixing $a/u>0$ and building a confidence interval for each sample $x_i$, proves -- not surprisingly -- the effectiveness of the $[x_i-u,x_i+u]$ interval. The test has been carried out by setting the constraints (\ref{constraint}-$b$); negative intervals have been allowed.

The generation of non-physical intervals lying entirely in the negative region is crucial to comply with the stipulated confidence level. As shown in  Fig.\ \ref{Neyman}, if these intervals are rejected and the measurements repeated until a physically acceptable interval is observed, the confidence level of the procedure is higher that what stated in (\ref{Paa}). Still worse, it is unpredictable, because it depends on the (unknown) measurand value.

\subsection{Credible intervals}
Since it must be conditional on a fixed measured value, a frequency-of-occurrence model of $\Prob(a\in [a_1,a_2]|x_0)$ must rely on the $\{a_i,x=\const\}$ sample space; that is, on the set of measurand value consistent with a unique measured value and credible interval. The sampling from this space can be carried out as follows. Firstly, a measurand value (say, $a_i$) is sampled according to $\pi(\phi)$ -- which encodes the pre-data information about $a$; secondly, a measurement results is sampled according to $P_x(\xi|a_i)$; thirdly, if the measurement result is $x_0$ -- in practice, to within some approximation -- the $a_i$ value is accepted; otherwise, it is rejected. In a long series of repetitions of this procedure, the fraction $\alpha$ of the accepted $a_i$ values is expected to be inside the single credible interval built by solving (\ref{CI3}).

\begin{figure}\centering
\includegraphics[width=7cm]{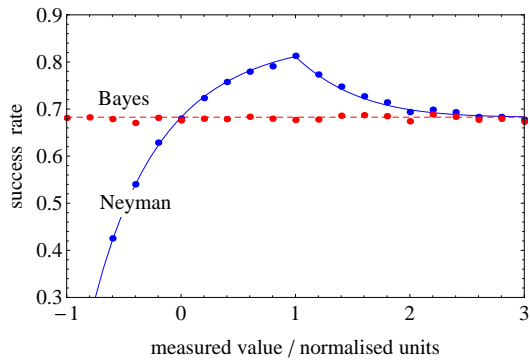}
\caption{Gaussian measurements of a positive quantity: a single confidence interval and a single credible interval have been built for the same measured value repeatedly sampled according to different distributions $P_x(\xi|a)$, each distribution corresponding to a different measurand value. Neyman: success rates of confidence intervals; Bayes: success rates of credible intervals. Lines are the theoretically expected values $\Prob(a\in[\au,\ao]|x_0)$ and $\Prob(a\in [a_1,a_2]|x_0)$; dots are the frequencies observed in Monte Carlo simulations. The calculations have been carried out by setting (\ref{constraint}-$b$) and (\ref{c-Bayes}-$b$).} \label{Bayes}
\end{figure}

In the case study here considered, the sampling from $\{a_i,x=\const\}$ and the assessment of the success rate of credible intervals can be carried out without assuming any prior distribution of the $a$ values by the following numerical experiment. Firstly, a measurand value (say, $a_0$) is chosen in whichever way; next, a measurement results (say, $x_i$) is sampled according to $P_x(\xi|a_0)$. To have $x_0$ instead, the measurand value is shifted to $a_i = a_0+x_0-x_i$. If $a_i<0$, the measurand value is rejected and the experiment is repeated. Otherwise, if $a_i>0$, it is checked if $a_i$ is in the fixed credible interval obtained by solving (\ref{CI3}). It is worth noting that a uniform prior distribution of the measurand values emerges naturally from the model, without being predetermined. The test has been carried out by setting the constraints (\ref{c-Bayes}-$b$); not even saying, as the horizontal line in Fig.\ \ref{Bayes} shows, the observed success rate is 0.68.

\subsection{Willink's paradox}
\begin{figure}\centering
\includegraphics[width=7cm]{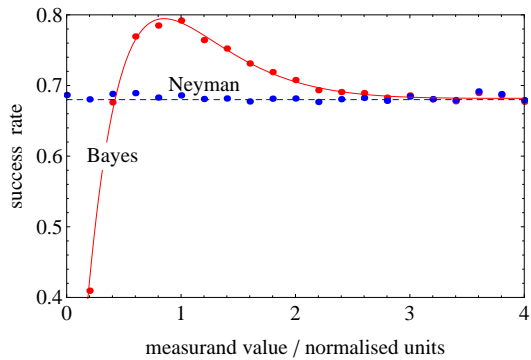}
\caption{Gaussian measurements of a positive quantity. Dots: success rates of confidence -- Neyman, the dashed line is the theoretically expected value $\Prob([\au,\ao]\ni a|a)$ -- and credible -- Bayes, the solid line is a smoothed interpolation of the data -- intervals. The calculations have been carried out by setting (\ref{constraint}-$b$) and (\ref{c-Bayes}-$b$), respectively. The sample space is $\{a=\const,x_i\}$: different credible- and confidence-intervals have been built from the results of repeated measurements of the same measurand.} \label{Neyman-2}
\end{figure}

The agrement of the long-run success rate of credible intervals with the predetermined coverage probability contradicts what observed by Willink \cite{Willink:2010a} (Fig.\ 4 -- solid line), which is reproduced by the solid line in Fig.\ \ref{Neyman-2}. This figure shows that, when $a/u \lesssim 1$, the success rate of the Bayesian intervals disagrees with the expected value. The paradox is solved by observing that, in Fig.\ \ref{Neyman-2} and \cite{Willink:2010a}, the success rate is calculated conditionally on the measurand value; that is, by fixing the measurand value and by building a new credible interval for each different measured value. This is equivalent to calculate $\Prob([a_1,a_2]\ni a|a)$. But Bayesian intervals are solutions of (\ref{Pax0}) and a frequency-of-occurrence model of $\Prob(a\in [a_1,a_2]|x_0)$ must be conditional on the measured value; that is, it must rely on a fixed measurement result. Therefore, the paradox originates from the use of the sample space $\{a=\const,x_i\}$.

To investigate further the differences between (\ref{Pax0}) and (\ref{Paa}), we calculated the probability -- $\Prob(a\in[\au,\ao]|x_0) = F_a(\ao|x_0)-F_a(\au|x_0)$, where $F_a(\phi|x_0)$ is given by (\ref{erf}) -- that the Neyman interval $[\au(x_0),\ao(x_0)]$ built from the measured value $x_0$ embeds the (unknown) value of the measurand. The result is shown in Fig.\ \ref{Bayes}. In addition, the success rate of $[\au,\ao]$ has been calculated by a Monte Carlo simulation conditional on a fixed measured value; that is, by sampling from the $\{a_i, x=\const\}$ space, where $x_i$, which is know, is fixed and $a$, which is unknown, is random. This has been done by the same numerical experiment used to assess the success rate of credible intervals; as expected, Fig.\ \ref{Bayes} shows the poor performance of the Neyman procedure when tested in this way.

\section{Conclusions}
This paper examined interval estimation from both the Neyman and Bayesian viewpoints and investigated differences not always perceived. It demonstrated a frequentist model of the coverage probability of Bayesian intervals, where a single interval is built for the same measured valued repeatedly sampled according to different distributions $P_x(\xi|a)$, each distribution corresponding to a different measurand value. No prior distribution of the measurand values has been explicitly assumed; rather, a uniform prior distribution emerges naturally from the model. Eventually, the paper proposed a solution to the paradoxical failure of the Bayesian intervals to pass a success-rate test based on repeated measurements of the same measurand and showed that an equivalent failure occurs when the Neyman intervals are tested against the same result repeatedly obtained by measuring different measurands.

The Neyman's view is the measurand value is not random; it is fixed and deterministic. Therefore, he discarded the specification of interval estimation given by (\ref{Pax0}) and turned to (\ref{Paa}). This attitude and the requirement that $\Prob([\au,\ao]\ni a|a)$ is independent of the $a$ value may lead to see the confidence level as the probability, $\Prob(a\in[\au,\ao]|x_0)$, of the measurand to be in a given $[\au,\ao]$ interval, rather than what it is, namely the probability, $\Prob([\au,\ao]\ni a|a)$, of a future $[\au,\ao]$ interval to encompass a given $a$ value.

Both the Neyman and Bayesian approaches are correct, but confidence and credible intervals are solutions of different problems, namely (\ref{Pax0}) and (\ref{Paa}). Hence, what is the best approach is an ill posed question. Whether to use one or the other to express the uncertainty of measurements depends on what problem we must solve and on decision theoretic considerations that are outside the scope of this paper. The following thoughts may supply some guidelines.

The construction of a generator of (random) intervals having a stipulated success rate of generating intervals including a fixed measurand, must rely on the Neyman procedure. For the Neyman's {\it practical statistician} in \cite{Neyman:1937}, the motivation of using confidence intervals lies in the customer satisfaction. If she sells confidence intervals, in the long run, she is sure that the fraction $\alpha$ of her customers had a correct statement.  But, the probability of the measurand to be in any specific interval may be not equal to the success-rate: in his seminal paper, Neyman already stressed that the confidence level is not the probability that the measurand is in the calculated interval.

If we want to express the measurement uncertainty by stating the probability that the measurand value is within a stipulated interval, which statement is not implied by Neyman intervals, we must rely on the probability distributions $\pi(\phi)$, before the measurement, and $P_a(\phi|x_0)$, after the measurement. The need of a prior distribution is an unavoidable consequence of the product rule of probabilities that discourages the use of credible intervals, because of lack of objectivity. However, an objectivity request does not make the Bayes theorem to vanish; to calculate the probability that the measurand is embedded in a given interval without the use of a prior distribution is impossible.

In order to allow the decision makers to make the relevant inferences by combining the result with any other information they have, it is incumbent on metrologists to provide the probability distribution $P_x(\xi|a)$ or, at least, the variance of the population of the possible results. But, if we must come to a decision based on the measurand value (e.g., to choose a value of the Planck constant to redefine the mass unit) it is the posterior probability density $P_a(\phi|x_0)$ -- hence, credible intervals -- that we need in order to maximise the expected utility (e.g., the continuity of the kilogram realizations). In addition, to account for the model uncertainty, we need also the evidence of the measurement result $Z(x_0)$.

\section*{Acknowledgements}
This work was jointly funded by the European Metrology Research Programme (EMRP) participating countries within the European Association of National Metrology Institutes (EURAMET) and the European Union.


\section*{References}
\begin{thebibliography}{99}
\bibitem{Neyman:1935}
 Neyman J 1935 On the problem of confidence intervals {\it Ann.\ Math.\ Stat.} {\bf 6} 111-6
\bibitem{Neyman:1937}
 Neyman J 1937 Outline of a theory of statistical estimation based on the classical theory of probability {\it Philos.\ Trans.\ Roy.\ Soc.\ Ser.\ A} {\bf 236} 333-80
\bibitem{Jaynes:1976}
 Jaynes E T 1976 Confidence intervals {\it vs.} Bayesian intervals in: {\it Foundations of Probability Theory, Statistical Inference, and Statistical Theories of Science} vol.\ II 175-257 (Dordrecht, Holland: D.\ Reidel Publishing Company)
\bibitem{Sivia}
 Sivia D S and Skilling J 2007 {\it Data Analysis: a Bayesian Tutorial} (Oxford: Oxford University Press)
\bibitem{Stein:1959}
 Stein C 1959 An example of a wide discrepancy between fiducial and confidence intervals {\it Ann.\ Math.\ Stat.} {\bf 30} 877-80
\bibitem{Feldman:1998}
 Feldman G J and Cousins R D 1998 Unified approach to the classical statistical analysis of small signals {\it Phys.\ Rev.\ D} {\bf 57} 3873-89
\bibitem{CERN:2000}
 D'Agostini G 2000 Confidence limits: What is the problem? Is there the solution {\it Workshop on confidence limits} eds.\ James F, Lyons L and Perrin Y (Genève: CERN)
\bibitem{Bukin:2003}
 Bukin A D 2003 A comparison of methods for confidence intervals {\it SLAC-R-703 Proceedings of PHYSTAT-2003: Statistical problems in particle physics, astrophysics and cosmology} eds.\ Lyons L, Mount R P and Reitmeyer R (Menlo Park: SLAC) 148-50
\bibitem{Lira:2006}
 Lira I and Woeger W 2006 Comparison between the conventional and Bayesian approaches to evaluate measurement data {\it Metrologia} {\bf 43} S249-59
\bibitem{Hall:2008}
 Hall B D 2008 Evaluating methods of calculating measurement uncertainty {\it Metrologia} {\bf 45} L5-8
\bibitem{Lira:2008}
 Lira I 2008 On the long-run success rate of coverage intervals {\it Metrologia} {\bf 45} L21-3
\bibitem{Wang-Iyer:2009}
 Wang C M and Iyer H K 2009 Fiducial intervals for the magnitude of a complex-valued quantity {\it Metrologia} {\bf 46} 81-6
\bibitem{Possolo:2009}
 Possolo A, Toman B and Estler T 2009 Contribution to a conversation about the Supplement 1 to the GUM {\it Metrologia} {\bf 46} L1-7
\bibitem{Willink:2010}
 Willink R 2010 On the validity of methods of uncertainty evaluation {\it Metrologia} {\bf 47} 80-9
\bibitem{Willink:2010c}
 Willink R 2010 Probability, belief and success rate: comments on 'On the meaning of coverage probabilities' {\it Metrologia} {\bf 47} 343-6
\bibitem{Attivissimo:2012}
 Attivissimo F, Giaquinto N and Savino M 2012 A Bayesian paradox and its impact on the GUM approach to uncertainty {\it Measurement} {\bf 45} 2194-202
\bibitem{Mana:2013}
 Bergamaschi L, D'Agostino G, Giordani L, Mana G and Oddone M 2013 The detection of signals hidden in noise {\it Metrologia} {\bf 50} 269-76
\bibitem{Calonico:2008}
 Calonico D, Levi F, Lorini L and Mana G 2009 Bayesian inference of a negative quantity from positive measurement results {\it Metrologia} {\bf 46} 267-71
\bibitem{Calonico:2009}
 Calonico D, Levi F, Lorini L and Mana G 2009 Bayesian estimate of the zero-density frequency of a Cs fountain {\it Metrologia} {\bf 46} 629-36
\bibitem{Willink:2010a}
 Willink R 2010 Uncertainty in repeated measurement of a small non-negative quantity: explanation and discussion of Bayesian methodology {\it Accred.\ Qual.\ Assur.} {\bf 15} 181-8
\bibitem{Willink:2010b}
 Willink R 2010 Measurement of small quantities: further observations on Bayesian methodology {\it Accred.\ Qual.\ Assur.} {\bf 15} 521-7
%\bibitem{Skilling:2006}
% Skilling J 2006 Nested sampling for general Bayesian computation {\it Bayesian Analysis} {\bf 1} 833-60
\end{thebibliography}
\end{document}